\documentclass[onepage,11pt,a4paper,pdftex]{article}

\usepackage{titletoc}
\usepackage{graphicx, amsfonts, amsmath, amssymb, amscd,  theorem, exscale,
epic, eepic, epsfig}

\usepackage{hyperref} %

\usepackage{makeidx}
\makeindex

\newcounter{Figure}
\theoremstyle{plain}
\theoremheaderfont{\scshape}

\newtheorem{Def}{\bf Definition}
\newtheorem{The}{\bf Theorem}

\theorembodyfont{\upshape}

\theoremheaderfont{\scshape\small} \theorembodyfont{\scshape\small}

\newcommand{\slap}{\mbox{$ \triangle \mkern -13mu / \ $}}
\newcommand{\nlap}{\mbox{$ \nabla \mkern -13mu / \ $}}
\newcommand{\dlap}{\mbox{$ div \mkern -13mu / \ $}}
\newcommand{\Dlap}{\mbox{$ D \mkern -13mu / \ $}}

\newcommand{\clap}{\mbox{$ curl \mkern -23mu / \ $}}

\newcommand{\be}{\begin{equation}}
\newcommand{\ee}{\end{equation}}
\newcommand{\bea}{\begin{eqnarray}}
\newcommand{\eea}{\end{eqnarray}}
\newcommand{\beas}{\begin{eqnarray*}}
\newcommand{\eeas}{\end{eqnarray*}}

\begin{document}

\begin{center}
{\Large \bf Null Asymptotics of Solutions of the Einstein-Maxwell Equations in
General Relativity 
\\ \vspace{5pt} 
and Gravitational Radiation} \\ 
\end{center}
\vspace{5pt}
\begin{center}
{\large \bf Lydia Bieri, PoNing Chen, Shing-Tung Yau } \\ 
\end{center}
\vspace{5pt}
\footnote{L. Bieri is supported by NSF grant DMS-0904583 
and S.-T. Yau is supported by NSF
grant PHY-0937443 and DMS-0904583. \\ 
Lydia Bieri, University of Michigan, Department of Mathematics, Ann Arbor MI. lbieri@umich.edu \\ 
PoNing Chen, Harvard University, Department of Mathematics, Cambridge MA. pchen@math.harvard.edu \\ 
Shing-Tung Yau, Harvard University, Department of Mathematics, Cambridge MA. yau@math.harvard.edu}
{\bf Abstract:} 
We prove that for spacetimes solving the Einstein-Maxwell (EM) equations, the
electromagnetic field 
contributes at highest order to the nonlinear memory effect of gravitational waves. 
In \cite{chrmemory} D. Christodoulou showed that gravitational waves have a
nonlinear memory. 
He discussed how this effect can be measured as a permanent displacement of test
masses in a laser interferometer gravitational wave detector. Christodoulou derived
a precise formula for this permanent displacement in the Einstein vacuum (EV) case. 
We prove in Theorem \ref{displ*1} that for the EM equations this permanent
displacement exhibits a term coming from the electromagnetic field. This term is at
the same highest order as the purely gravitational term that governs the EV
situation. 
On the other hand, in Chapter \ref{wave}, we show that to leading order, the
presence of the electromagnetic field does not change the instantaneous displacement
of the test masses. 
Following the method introduced by D. Christodoulou in \cite{chrmemory} and
asymptotics derived by N. Zipser in \cite{zip} and \cite{zip2}, 
we investigate gravitational radiation at null infinity in spacetimes solving the EM
equations. 
We study the Bondi mass loss formula at null infinity derived in \cite{zip2}. 
We show that 
the mass loss formula from \cite{zip2}  is compatible with the one in Bondi
coordinates obtained in \cite{vanderBurg}. 
And we observe that the presence of the electromagnetic field increases the total
energy radiated to infinity up to leading order. 
Moreover, we compute the limit of the area radius at null infinity in Theorem
\ref{r*1}. 
\tableofcontents

\section{Introduction and Main Results}
In this paper we investigate the null asymptotics for spacetimes solving the
Einstein-Maxwell (EM) equations, 
compute the radiated energy, derive limits at null infinity and compare them with
the Einstein vacuum (EV) case. We show that the presence of the 
electromagnetic field does not affect the instantaneous displacement of the test
masses of a laser interferometer detector at 
leading order, as it only comes in at lower order. But the electromagnetic field
does contribute to the nonlinear effect of 
the displacement of the test masses. 
The EM case gives us a wonderful opportunity to 
observe mass loss and also to measure gravitational radiation. It is crucial to
understand fully the 
behavior of the gravitational field also when other fields are present and to
investigate their interplay. 
The only way to achieve this, is to compute the null asymptotics of the
spacetimes.\\ \\
A major goal of mathematical General Relativity (GR) and astrophysics is to
precisely describe and 
finally observe gravitational radiation, one of the predictions of GR. In order to
do so, one has to study 
the null asymptotical limits of the spacetimes for typical sources. 
Among the latter we find binary neutron stars and binary black hole mergers. 
In these processes typically mass and momenta are radiated away in form of
gravitational waves. 
Bondi, van der Burg and Metzner studied these in \cite{BBM}. 
D. Christodoulou showed in his paper \cite{chrmemory} that every gravitational-wave
burst has a 
nonlinear memory. The insights of this work are based on the precise description of
null infinity obtained 
by D. Christodoulou and S. Klainerman in \cite{sta}. Among the many pioneering
results they derived the 
Bondi mass loss formula. This is all in the regime of the Einstein vacuum equations. 
Then N. Zipser studied the Einstein-Maxwell equations in \cite{zip}, \cite{zip2} and
computed limits 
along the lines of \cite{sta} for this case. She derived a Bondi mass loss formula,
where in addition to the one 
obtained by Christodoulou and Klainerman, a component of the electromagnetic field
comes in. 
Thus the mass radiated away goes into the gravitational and the electromagnetic field. 
Here, we rely on the methods introduced 
in \cite{sta}, used in \cite{zip}, \cite{zip2} and by one of the present authors in
\cite{lydia1}, \cite{lydia2}. 
There is a large literature about gravitational radiation. However, in the present
paper, we only give the references which are 
relevant to our investigations.  \\ \\ 
The main results of this paper are the following. We first recall the Bondi mass
loss formula obtained in \cite{zip2} for spacetimes solving the EM equations. 
\begin{equation*}
\frac{\partial }{\partial u}M\left( u\right) =\frac{1}{8\pi }%
\int_{S^{2}}\left( \left| \Xi \right| ^{2}+\frac{1}{2}\left| A_{F}\right|
^{2}\right) d\mu _{\overset{\circ }{\gamma }}
\end{equation*}
Compared to the formula obtained in \cite{sta} for spacetimes solving the EV
equations, we have an additional term, $|A_F|^2$, from the electromagnetic field.
Furthermore,
we compare the above mass loss formula to the corresponding formula in Bondi
coordinates \cite{vanderBurg} 
$$
\frac{\partial }{\partial w}M(w) = - \int_{S^2} \left( (\partial_w c)^2 +(\partial_w
d)^2  +\frac{1}{2}(X^2+Y^2)  \right) d\mu _{\overset{\circ }{\gamma }}
$$ 
and show that the two formulae agree.  
 \\  \\
As shown in the work of Christodoulou \cite{chrmemory}, $\Sigma^+ - \Sigma^-$ is the
term which governs the permanent displacement of test particles. Using this fact, 
Christodoulou shows that the gravitational field has a non-linear ``memory" which
can be detected by a gravitational-wave experiment in a spacetime solving the EV
equations. 
We will describe this experiment in the last section as well. In section
\ref{permanentdispl***} of our paper, we study the permanent displacement formula
for uncharged test particles of the same gravitational-wave experiment in a
spacetime solving the EM equations and show that the electromagnetic field
contributes to the nonlinear effect. 
We first obtain Theorem \ref{displ*1} which determines $\Sigma^+ - \Sigma^-$ in the
EM case. From Theorem \ref{displ*1}, we observe that the electromagnetic field
changes the leading order term of the permanent displacement of test particles. Then
in the last section, we study in details a gravitational wave experiment for our
findings.  We observe that the electromagnetic field  does not enter the leading
order term of the Jacobi equation. As a result, to leading order, it does not change
the instantaneous displacement of test particles. 
But the electromagnetic field does contribute at highest order to the nonlinear
effect of the permanent displacement of test masses. 
Furthermore, in Theorem \ref{r*1}  we compute the limit of the area radius $r$ on
any null hypersurface $C_u$ as $t$ goes to $\infty$ and show that the result
coincides with the one obtained in \cite{sta} for EV.
\\  \\ 
We follow the method  introduced  by Christodoulou in \cite{chrmemory} to study the
effect of gravitational waves. The treatment is based on the asymptotic behavior of
the gravitational field obtained at null and spatial infinity. These rigorous
asymptotics allow us to study the structure of the spacetimes at null infinity. 
The spacetime is foliated by a time function $t$ and by an optical function $u$. 
The corresponding lapse functions are denoted by $\phi$ and $a$. 
Each level set of $u$, $C_u$, is an outgoing null hypersurface and each level set of
$t$, $H_t$ is a maximal spacelike hypersurface.  We pick a suitable pair of normal
vectors along the null hypersurface. The flow along these vector fields generates a
family of diffeomorphisms $\phi_u$ of $S^2$. We use $\phi_u$ to pull back tensor
fields in our spacetime. This allows us to study their limit at null infinity along
the null hypersurface $C_u$.  Then we study the effect of gravitational waves by
taking the limit as $u$ goes to $\pm \infty$. Christodoulou in \cite{chrmemory}
gives a complete explanation of the structure at null infinity.
\\ \\
The methods introduced in \cite{sta}, used in \cite{zip}, \cite{zip2}, reveal the
structure of 
the null asymptotics of our spacetimes. 
In these works, to prove the stability result, the data is assumed to be small.
However, as far as the 
study of the null asymptotics is concerned, the data can be large. 
We give a brief outline in the last part of this introduction 
of the methods of \cite{sta}. \\ \\
In General Relativity the fundamental equations are the Einstein equations 
linking the curvature of the spacetime to its matter content. 
\be \label{equsGRgRT}
G_{\mu \nu} \  := \ R_{\mu \nu} \ - \ \frac{1}{2} \ g_{\mu \nu} \ R \ = \ 8 \pi \
T_{\mu \nu}  \ , 
\ee
(setting $G = c  = 1$), $\mu, \nu = 0,1,2,3$.  
$G_{\mu \nu}$ denotes the Einstein tensor,  
$R_{\mu \nu}$ is the Ricci curvature tensor,  
$R$ the scalar curvature tensor, 
$g$ the metric tensor and 
$T_{\mu \nu}$ is the energy-momentum tensor. \\ \\ 
Here, we are discussing the Einstein-Maxwell equations. 
This means that $T_{\mu \nu}$ on the right hand side of (\ref{equsGRgRT}) is the 
stress-energy tensor of the electromagnetic field. 
The twice contracted Bianchi identities imply that 
\be \label{DG01}
D^{\nu} G_{\mu \nu} = 0 \ \ . 
\ee
This is equivalent to the following equation, namely, that the divergence of the
stress-energy tensor of the electromagnetic 
field vanishes: 
\be \label{DT01}
D^{\nu} T_{\mu \nu} = 0 
\ee
with 
\be \label{T*1}
T_{\mu \nu} = \frac{1}{4 \pi} \big(  F_{\mu}^{\ \rho} F_{\nu \rho} - \frac{1}{4}
g_{\mu \nu} F_{\rho \sigma} F^{\rho \sigma}  \big)
\ee
where $F$ denotes the electromagnetic field. 
Note that $F$ is an antisymmetric covariant 2-tensor. 
As $T_{\mu \nu}$ is trace-free, the Einstein equations (\ref{equsGRgRT}) take the form 
\be \label{equsGRgEM}
R_{\mu \nu} = 8 \pi T_{\mu \nu} \ \ . 
\ee
We find that the scalar curvature is identically zero. We write 
the Einstein-Maxwell (EM) equations as 
\bea
R_{\mu \nu} & =  & 8 \pi T_{\mu \nu}  \\ 
D^{\alpha} F_{\alpha \beta} & = & 0 \\ 
D^{\alpha} \ ^* F_{\alpha \beta} & = & 0. \ \ 
\eea
As a consequence of the Maxwell equations, we have
\begin{equation}
\Box F =0
\end{equation}
where $\Box$ is the de Rham Laplacian with respect to the metric $g$. \\ \\
We split the Riemannian curvature $R_{\alpha \beta \gamma \delta}$ into its
traceless part, namely the 
Weyl tensor $W_{\alpha \beta \gamma \delta }$, and a part including the spacetime
Ricci curvature $R_{\alpha \beta}$ 
and spacetime scalar curvature $R$: 
\bea
R_{\alpha \beta \gamma \delta } &=&W_{\alpha \beta \gamma
\delta }+\frac{1}{2}(g_{\alpha \gamma }R_{\beta \delta }+
g_{\beta \delta }R_{\alpha \gamma }-g_{\beta
\gamma }R_{\alpha \delta }-g_{\alpha \delta }R_{\beta \gamma })  \nonumber \\
&&-\frac{1}{6}\left( g_{\alpha \gamma }g_{\beta \delta }- 
g_{\alpha \delta }g_{\beta \gamma }\right) R  \ \ . \label{Riem*tracelessRicci}
\eea
One can then define the Bel-Robinson tensor 
\be
Q_{\alpha \beta \gamma \delta} = W_{\alpha \rho \gamma \sigma} W_{\beta \ \delta}^{\
\rho \  \sigma} + 
\ ^* W_{\alpha \rho \gamma \sigma} \ ^* W_{\beta \ \delta}^{\ \rho \  \sigma}  \ \ . 
\ee
The Bianchi equations for the Weyl tensor in the presence of an electromagnetic
field read 
\be
D^{\alpha} W_{\alpha \beta \gamma \delta} = \frac{1}{2} ( D_{\gamma} R_{\beta
\delta} - D_{\delta} R_{\beta \gamma} ) \ \ . 
\ee
In \cite{sta} Christodoulou and Klainerman derived the asymptotic behavior in the
case of 
strongly asymptotically flat initial data of the following type. 
\begin{Def} \label{CKID} 
We define a strongly asymptotically flat initial data set in the sense 
of \cite{sta} (studied by Christodoulou and Klainerman) to be 
an initial data set $(H, \bar{g}, k)$, where 
$\bar{g}$ and $k$ are sufficiently smooth and there exists a coordinate system 
$(x^1, x^2, x^3)$ defined in a neighbourhood of infinity such that, 
as  $r = (\sum_{i=1}^3 (x^i)^2 )^{\frac{1}{2}} \to \infty$,  $\bar{g}_{ij}$ and
$k_{ij}$ are: 
\bea
\bar{g}_{ij} \ & = & \ (1 \ + \ \frac{2M}{r}) \ \delta_{ij} \ + \ o_4 \ (r^{-
\frac{3}{2}}) \label{safg} \\ 
k_{ij} \ & = & \  o_3 \ (r^{- \frac{5}{2}}) \ ,  \label{safk}
\eea
with $M$ denoting the mass. 
\end{Def}
Under a smallness condition on the initial data, Christodoulou and Klainerman proved
in \cite{sta} that 
this can be extended uniquely to a smooth, globally hyperbolic and geodesically
complete spacetime solving the EV equations. 
The resulting spacetime is globally asymptotically flat. Together with the existence
and uniqueness theorem comes 
a precise description of the asymptotic behavior of the spacetime. 
While the smallness condition was imposed in order to ensure completeness, the
results about the 
behavior at null infinity are largely independent of the smallness. 
The decay behavior of the components of the Weyl tensor are given below. And the
limits at null infinity of the 
relevant quantities are given in section \ref{chnullasymptotics}. \\ \\ 
In \cite{sta} as well as in \cite{zip}, \cite{zip2} and \cite{lydia1}, \cite{lydia2}
the Weyl tensor $W$ in $(M,g)$ 
is decomposed with respect to the 
null frame $e_4, e_3, e_2, e_1$. That is, $e_4$ and $e_3$ form a 
null pair which is supplemented by $e_A, \ A = 1, 2$, a local 
frame field for $S_{t,u} = H_t \cap C_u $. Given this null pair, $e_3$ and $e_4$, we
can define the tensor 
of projection from the tangent space of $M$ to that of $S_{t,u}$.
\[  \Pi^{\mu \nu} = g^{\mu \nu} + 
\frac{1}{2}(e_4^{\nu}e_3^{\mu}+e_3^{\nu}e_4^{\mu}). \]
We decompose the second fundamental form $k_{ij}$ of $H_t$ into
\bea
k_{NN} &=& \delta \\
k_{AN} &=& \epsilon_A\\
k_{AB} &=& \eta_{AB} 
\eea
where $N$ is the unit normal vector of $S_{t,u}$ in $H_t$. Let $T$ be the future-directed unit normal to $H_t$.
We define
\[  \theta_{AB} =  \langle \nabla_{A} N, e_B \rangle. \]
The Ricci coefficients of the null standard frame $T-N,T+N,e_2, e_1$ are given by
the following
\bea
\chi^{\prime}_{AB} &=& \theta_{AB} - \eta_{AB} \\
\underline{\chi}^{\prime}_{AB} &=& - \theta_{AB} - \eta_{AB} \\
\underline{\xi}^{\prime}_A &= &\phi^{-1} \nlap_A \phi - a^{-1} \nlap_A a \\
\underline{\zeta}^{\prime}_A &= &\phi^{-1} \nlap_A \phi - \epsilon_A  \\
\zeta^{\prime}_A &= &\phi^{-1} \nlap_A \phi + \epsilon_A  \\
\nu^{\prime} &= & - \phi^{-1} \nlap_N \phi + \delta  \\
\underline{\nu}^{\prime} &= &  \phi^{-1} \nlap_N \phi + \delta   \\
\omega^{\prime} &=&\delta -a^{-1} \nlap_N a
\eea
We use $\chi$, $\underline{\chi}$, etc for the  Ricci coefficients of the null frame
$a^{-1}(T-N),a(T+N),e_2, e_1$.
\begin{Def}
We define the null components of $W$ as follows: 
\bea
\underline{\alpha}_{\mu \nu} \ (W) \ & = & \ 
\Pi_{\mu}^{\ \rho} \ \Pi_{\nu}^{\ \sigma} \ W_{\rho \gamma \sigma \delta} \
e_3^{\gamma} \ e_3^{\delta} 
\label{underlinealpha} \\ 
\underline{\beta}_{\mu} \ (W) \ & = & \ 
\frac{1}{2} \ \Pi_{\mu}^{\ \rho} \ W_{\rho \sigma \gamma \delta} \  e_3^{\sigma} \
e_3^{\gamma} \ e_4^{\delta} 
\label{underlinebeta} \\ 
\rho \ (W) \ & = & \ 
\frac{1}{4} \ W_{\alpha \beta \gamma \delta} \ e_3^{\alpha} \ e_4^{\beta} \
e_3^{\gamma} \ e_4^{\delta} 
\label{rho} \\ 
\sigma \ (W) \ & = & \ 
\frac{1}{4} \ \ ^*W_{\alpha \beta \gamma \delta} \ e_3^{\alpha} \ e_4^{\beta} \
e_3^{\gamma} \ e_4^{\delta} 
\label{sigma} \\ 
\beta_{\mu}  \ (W) \ & = & \  
\frac{1}{2} \ \Pi_{\mu}^{\ \rho} \ W_{\rho \sigma \gamma \delta} \ e_4^{\sigma} \
e_3^{\gamma} \ e_4^{\delta} 
\label{beta} \\ 
\alpha_{\mu \nu} \ (W) \ & = & \ 
\Pi_{\mu}^{\ \rho} \ \Pi_{\nu}^{\ \sigma} \ W_{\rho \gamma \sigma \delta} \
e_4^{\gamma} \ e_4^{\delta}  \ . 
\label{alphaR}
\eea
\end{Def}
The estimates in \cite{sta} yield the decay behavior: 
\beas 
\underline{\alpha}(W) \ & = & \ O \ (r^{- 1} \ \tau_-^{- \frac{5}{2}}) \\ 
\underline{\beta}(W) \ & = & \ O \ (r^{- 2} \ \tau_-^{- \frac{3}{2}}) \\ 
\rho(W) \   & = & \ O \ (r^{-3})  \\ 
\sigma(W)  \ & = & \ O \ (r^{-3} \ \tau_-^{- \frac{1}{2}}) \\ 
\alpha(W) , \ \beta(W) \ & = & \ o \ ( r^{- \frac{7}{2}})  
\eeas
where $\tau_-^2=1+u^2$ and $r(t,u)$ is the area radius of the surface $S_{t,u}$. \\ \\
In  \cite{zip}, \cite{zip2}, Zipser works with the same conditions on the metric,
second fundamental form and curvature, 
in addition she imposes a decay
condition on the electromagnetic field $F$, namely 
\begin{equation}
F|_{H }=o_{3}\left( r^{-\frac{5}{2}}\right) .  \label{initial F}
\end{equation}
The null components\ of the electromagnetic field are written as 
\begin{equation}
\begin{array}{lll}
F_{A3}=\underline{\alpha }(F)_{A} &  & F_{A4}=\alpha \left( F\right) _{A} \\ 
F_{34}=2\rho \left( F\right) &  & F_{12}=\sigma \left( F\right) .%
\end{array}
\label{null-electricfield}
\end{equation}%
The corresponding null decomposition $\left \{ \underline{\alpha }%
\left( ^{\ast }F\right) ,\alpha \left(^{\ast }F\right) ,\rho \left(
^{\ast }F\right) ,\sigma \left(  ^{\ast }F \right) \right \} $ of $^{\ast }F$ is
given by 
\begin{eqnarray}
\underline{\alpha }\left( ^{\ast }F\right) _{A} &=&-\underline{\alpha }%
\left( F\right) ^{B}\epsilon _{BA}\quad \quad \alpha \left(
^{\ast }F \right) _{A}=\alpha \left( F\right) ^{B}\epsilon _{BA}  \notag \\
\rho \left(   ^{\ast }F  \right) &=&\sigma \left( F\right) \quad \quad \, \,
\, \, \, \quad \quad \sigma \left(  ^{\ast }F  \right) =-\rho \left( F\right)
\end{eqnarray}%
where the Hodge dual of a tensor $u$ tangent to $S_{t,u}$, is defined by 
\begin{equation*}
^{\ast }u_{A}={\epsilon _{A}}^{B}u_{B.}
\end{equation*}%
The estimates in \cite{zip}, \cite{zip2} yield the decay behavior: 
\beas 
\underline{\alpha}(F) \ & = & \ O \ (r^{- 1} \ \tau_-^{- \frac{3}{2}}) \\ 
\rho(F), \sigma(F)  \   & = & \ O \ (r^{-2} \ \tau_-^{- \frac{1}{2}} )  \\ 
\alpha(F) \ & = & \ o \ ( r^{- \frac{5}{2}}).  
\eeas
One of the main difficulties in \cite{sta} is that a general spacetime has no
symmetries and thus does not have suitable vectorfields to construct integral 
conserved quantities. To overcome this difficulty, Christodoulou and Klainerman 
use the `closeness' of their spacetimes to the Minkowski spacetime and construct 
quasi-conformal vector fields. The main step is carried out within a bootstrap
argument in the 
`last slice', namely in a spacelike hypersurface which is a level set of the time
function $t$. 
First, the authors foliate the spacetime by functions $t$ and $u$ near the initial
slice.  From the 
foliations, one constructs vectorfields that are almost Killing. Combining these
vectorfields with the Bel-Robinson tensor, one obtains local estimates for 
the Weyl curvature tensor $W$ and the  electromagnetic field $F$. With these
estimates, one constructs a new optical function  which is defined on a 
larger domain in the spacetime.  Then, following a continuity argument, one obtains
a smooth, globally hyperbolic and geodesically complete spacetime
 solving the Einstein equations. The resulting spacetime is globally asymptotically
flat, satisfying the above decay properties. \\ \\ 
{\it \bf Acknowledgment.} We thank Demetrios Christodoulou for helpful discussions, critcal reading of our
paper and constructive feedback. 
\section{Null Asymptotics}
\label{chnullasymptotics}
\subsection{Asymptotic Behavior and Bondi Mass}
We need precise data at null infinity. 
In particular, we have to know the Bondi mass and the asymptotic behavior of the
components of the curvature and the electromagnetic field. 
Zipser described them in \cite{zip}, \cite{zip2} 
following the
discussion in Chapter 17 of \cite{sta}, making changes as necessary due to the
presence of
the electromagnetic field. The parameters of the foliations and the
components of the Weyl tensor behave exactly as in \cite{sta}. \ That
is, 
the following holds:  
along the null hypersurfaces $C_{u}$ as $t\rightarrow \infty $, it is  
\begin{equation}
\lim_{C_{u},t\rightarrow \infty }\phi =1,\, \, \, \ \ \ \ \ \ \ \ \ \ \ \,
\, \ \ \ \, \, \ \, \, \, \ \ \ \ \, \ \, \lim_{C_{u},t\rightarrow \infty
}a=1  \label{aymptotic behvouir of phi and a}
\end{equation}%
and 
\begin{equation}
\lim_{C_{u},t\rightarrow \infty }\left( rtr\chi \right) =2,\, \, \ \ \, \ \
\ \ \ \ \ \ \ \ \, \ \ \, \, \ \ \, \ \, \, \lim_{C_{u},t\rightarrow \infty
}\left( rtr\underline{\chi }\right) =-2  \label{asymptotic behaviour of trxi}
\end{equation}
Furthermore, we let 
\begin{equation}
H=\lim_{C_{u},t\rightarrow \infty }\left( r^2(tr\chi' -\frac{2}{r}) \right ) \, .
\end{equation}
From the existence theorem of \cite{zip}, \cite{zip2}, Zipser makes the following
conclusions, which are
generalizations of conclusions 17.0.1 through 17.0.4 in \cite{sta}. \\ \\ 
Following the convention in \cite{sta} and \cite{zip2}, the
pointwise norms $\left \vert \quad \right \vert $ of the tensors on $S^{2}$
relate to the metric $\overset{\circ }{\gamma }$, which is the limit of the
induced metrics on $S_{t,u}$ rescaled by $r^{-2}$ for each $u$ as $t\rightarrow \infty $ . 
\begin{The} \label{conclcurvandfieldcompts1}
On any null hypersurface $C_{u}$, the normalized curvature components $r%
\underline{\alpha }\left( W\right) $, $r^{2}\underline{\beta }\left(
W\right) $, $r^{3}\rho \left( W\right) $, $r^{3}\sigma \left( W\right) $, $r%
\underline{\alpha }\left( F\right) $, $r^{2}\rho \left( F\right) $, $%
r^{2}\sigma \left( F\right) $ have limits as $t\rightarrow \infty $, in particular 
\begin{eqnarray*}
\lim_{C_{u},t\rightarrow \infty }r\underline{\alpha }\left( W\right)
&=&A_{W}\left( u,\cdot \right) ,\, \ \ \ \ \ \ \ \ \ \ \ \
\lim_{C_{u},t\rightarrow \infty }\,r^{2}\underline{\beta }\left( W\right)
=B_{W}\left( u,\cdot \right) \\
\lim_{C_{u},t\rightarrow \infty }r^{3}\rho \left( W\right) &=&P_{W}\left(
u,\cdot \right) ,\, \ \ \ \ \ \ \ \ \ \ \ \ \lim_{C_{u},t\rightarrow \infty
}r^{3}\sigma \left( W\right) =Q_{W}\left( u,\cdot \right) \\
\lim_{C_{u},t\rightarrow \infty }r\underline{\alpha }\left( F \right)
&=&A_{F}\left( u,\cdot \right) , \\
\lim_{C_{u},t\rightarrow \infty }r^{2}\rho \left( F\right) &=&P_{F}\left(
u,\cdot \right) ,\, \ \ \ \ \ \ \ \ \ \ \ \ \lim_{C_{u},t\rightarrow \infty
}r^{2}\sigma \left( F\right) =Q_{F}\left( u,\cdot \right)
\end{eqnarray*}
with $A_{W}$ a symmetric traceless covariant 2-tensor, $B_{W}$ and $%
A_{F} $ 1-forms and $P_{W}$, $Q_{W}$, $P_{F}$, $Q_{F}$ functions on $%
S^{2}$ depending on $u$. The following decay properties hold: 
\begin{eqnarray*}
\left| A_{W}\left( u,\cdot \right) \right| &\leq &C\left( 1+\left| u\right|
\right) ^{-5/2}\, \, \ \ \ \ \ \ \ \ \ \ \ \ \left| B_{W}\left( u,\cdot
\right) \right| \leq C\left( 1+\left| u\right| \right) ^{-3/2} \\
\left| P_{W}\left( u,\cdot \right) -\overline{P}_{W}\left( u\right) \right|
&\leq &\left( 1+\left| u\right| \right) ^{-1/2}\, \, \quad \quad \quad \quad
\, \left| Q_{W}\left( u,\cdot \right) -\overline{Q}_{W}\left( u\right)
\right| \leq \left( 1+\left| u\right| \right) ^{-1/2} \\
\left| A_{F}\left( u,\cdot \right) \right| &\leq &C\left( 1+\left| u\right|
\right) ^{-3/2} \\
\left| P_{F}\left( u,\cdot \right) \right| &\leq &\left( 1+\left| u\right|
\right) ^{-1/2}\, \, \ \ \ \ \ \ \ \ \ \ \ \ \, \, \, \, \left| Q_{F}\left(
u,\cdot \right) \right| \leq \left( 1+\left| u\right| \right) ^{-1/2}
\end{eqnarray*}
and 
\begin{equation*}
\lim_{u\rightarrow -\infty }\overline{P}_{W}\left( u\right) =0,\, \ \ \ \ \
\ \lim_{u\rightarrow -\infty }\overline{Q}_{W}\left( u\right) =0.
\end{equation*}
\end{The}
The existence of the limits in the conclusion follows from the estimates in
the existence theorem of \cite{zip2} (i.e. \cite{zip}). 
\begin{The} \label{conclSigma1}
On the null hypersurface $C_{u}$, the normalized shear $r^{2} \widehat \chi ^{\prime }$
has limit as $t\rightarrow \infty $: 
\begin{equation*}
\lim_{C_{u},t\rightarrow \infty }r^{2}  \widehat \chi ^{\prime }=\Sigma \left( u,\cdot
\right)
\end{equation*}
with $\Sigma $ being a symmetric traceless covariant 2-tensor on $S^{2}$
depending on $u$.
\end{The}
The proof is the same as in \cite{sta} because the 
propagation equation stays unaltered 
\begin{equation*}
\frac{d\widehat{\chi }_{AB}}{ds}=-tr\chi \widehat{\chi }_{AB}-\alpha
(W)_{AB}.
\end{equation*}
\begin{The} \label{conclXi1}
On any null hypersurface $C_{u}$, the limit of $r\widehat{\eta }$ exists as $%
t\rightarrow \infty $, that is 
\begin{equation*}
\lim_{C_{u},t\rightarrow \infty }r\widehat{\eta }=\Xi \left( u,\cdot \right)
\end{equation*}
with $\Xi $ being a symmetric traceless 2-covariant tensor on $S^{2}$
depending on $u$ and having the decay property 
\begin{equation*}
\left| \Xi \left( u,\cdot \right) \right| _{\overset{\circ }{\gamma }}\leq
C\left( 1+\left| u\right| \right) ^{-3/2}.
\end{equation*}
Further, it is  
\begin{equation*}
\lim_{C_{u},t\rightarrow \infty }r\widehat{\theta }=-\frac{1}{2}%
\lim_{C_{u},t\rightarrow \infty }r\widehat{\underline{\chi }}^{\prime }=\Xi
\end{equation*}
as well as 
\begin{eqnarray}
\frac{\partial \Sigma }{\partial u} &=&-\Xi   \label{Sigmau*1}  \\
\frac{\partial \Xi }{\partial u} &=&-\frac{1}{4}A_{W}.  \label{Xiu*1}
\end{eqnarray}
\end{The}
Zipser proves this result as conclusion 3 in \cite{zip2}. The argument 
is along the lines of the proof of conclusion 17.0.3 in \cite{sta}. \\ \\ 
Zipser follows \cite{sta} to derive the Bondi mass formula by calculating a
propagation equation for the Hawking mass enclosed by a 2-surface $S_{t,u}$.
The Hawking mass is defined as 
\begin{equation}
m\left( t,u\right) =\frac{r}{2}\left( 1+\frac{1}{16\pi }\int_{S_{t,u}}tr\chi
tr\underline{\chi }\right) .  \label{hawking mass}
\end{equation}%
Let 
\begin{equation}
\underline{\mu }=-\dlap \underline{\zeta }+\frac{1}{2}\widehat{\chi }%
\cdot \underline{\widehat{\chi }}-\rho \left( W\right) -\frac{1}{2}\left(
\rho ^{2}\left( F\right) +\sigma ^{2}\left( F\right) \right) .
\label{Mass aspect function}
\end{equation}%
With respect to the $l$-pair, one has the null structure equations 
\begin{eqnarray*}
\frac{dtr\underline{\chi }}{ds}+\frac{1}{2}tr\chi tr\underline{\chi } &=&-2%
\underline{\mu }+2\left \vert \zeta \right \vert ^{2} \\
\frac{dtr\chi }{ds}+\frac{1}{2}\left( tr\chi \right) ^{2} &=&-\left \vert 
\widehat{\chi }\right \vert ^{2}-\left \vert \alpha \left( F\right) \right
\vert ^{2}.
\end{eqnarray*}%
One computes,  
\begin{eqnarray*}
\frac{d}{ds}tr\chi tr\underline{\chi }+tr\chi \left( tr\chi tr\underline{%
\chi }\right) &=&-2\underline{\mu }tr\chi +2tr\chi \left \vert \zeta \right
\vert ^{2} \\
&&-tr\underline{\chi }\left \vert \widehat{\chi }\right \vert ^{2}-tr%
\underline{\chi }\left \vert \alpha \left( F\right) \right \vert ^{2}
\end{eqnarray*}%
thus  
\begin{eqnarray}
\frac{\partial }{\partial t}\int_{S_{t,u}}tr\chi tr\underline{\chi }
&=&-2\int_{S_{t,u}}a\phi \underline{\mu }tr\chi
\label{t-integrate trxtrxbar} \\
&&+\int_{S_{t,u}}a\phi \left( -tr\underline{\chi }\left \vert \widehat{\chi }%
\right \vert ^{2}-tr\underline{\chi }\left \vert \alpha \left( F\right)
\right \vert ^{2}+2tr\chi \left \vert \zeta \right \vert
^{2}\right) .  \notag
\end{eqnarray}%
Using the Gauss equation 
\begin{equation*}
K=-\frac{1}{4}tr\chi tr\underline{\chi }+\frac{1}{2}\widehat{\chi }\cdot 
\underline{\widehat{\chi }}-\rho \left( W\right) -\frac{1}{2}\left( \rho
^{2}\left( F\right) +\sigma ^{2}\left( F\right) \right) ,
\end{equation*}%
one derives 
\begin{equation}
\underline{\mu }=- \dlap \underline{\zeta }+K+\frac{1}{4}tr\chi tr%
\underline{\chi }.  \label{Gauss mass aspect}
\end{equation}%
By the Gauss-Bonnet formula and formulas (\ref{Mass aspect function}), 
(\ref{Gauss mass aspect}), conclude that 
\begin{eqnarray}
\int_{S_{t,u}}\underline{\mu } &=&\int_{S_{t,u}}\left( \frac{1}{2}\widehat{%
\chi }\cdot \underline{\widehat{\chi }}-\rho \left( W\right) -\frac{1}{2}%
\left( \rho ^{2}\left( F\right) +\sigma ^{2}\left( F\right) \right) \right)
\label{integrate mu} \\
&=&4\pi \left( 1+\frac{1}{16\pi }\int_{S_{t,u}}tr\chi tr\underline{\chi }%
\right) =\frac{8\pi }{r}m.  \notag
\end{eqnarray}%
Moreover, by 
\begin{equation*}
\frac{d}{dt}r=\frac{r}{2}\overline{\phi atr\chi },
\end{equation*}%
and (\ref{t-integrate trxtrxbar}), (\ref{integrate mu}), it is  
\begin{eqnarray}
\frac{\partial }{\partial t}m\left( t,u\right) &=&-\frac{r}{16\pi }%
\int_{S_{t,u}}\left( a\phi tr\chi -\overline{\phi atr\chi }\right) 
\underline{\mu }  \label{dt-hawking mass} \\
&&+\frac{r}{8\pi }\int_{S_{t,u}}a\phi \left( \frac{1}{2}tr\chi \left \vert
\zeta \right \vert ^{2}-\frac{1}{4}tr\underline{\chi }\left \vert \widehat{%
\chi }\right \vert ^{2}-\frac{1}{4}tr\underline{\chi }\left \vert \alpha
\left( F\right) \right \vert ^{2}\right) .  \notag
\end{eqnarray}%
Note that $K+\frac{1}{4}tr\chi tr\underline{\chi }=O\left( r^{-3}\right) $, $%
\underline{\mu }=O\left( r^{-3}\right) $. From the asymptotic behavior of
the right-hand side of (\ref{dt-hawking mass}), it follows 
\begin{equation*}
\frac{\partial }{\partial t}m\left( t,u\right) =O\left( r^{-2}\right) .
\end{equation*}%
This means that $m\left( t,u\right) $ has a limit for any fixed $u$ as $%
t\rightarrow \infty $, namely the Bondi mass of the null hypersurface $C_{u}$. 
As in \cite{sta}, it is denoted by $M\left(u\right) $. 
The terms appearing due to the presence of the electromagnetic field are shown to
decay fast enough so that the mass decays at the same rate as in 
\cite{sta}. \ In particular, 
\begin{equation*}
m\left( t,u\right) =M\left( u\right) +O\left( r^{-1}\right)
\end{equation*}%
as $t\rightarrow \infty $ on $C_{u}$. \\ \\ 
Following \cite{sta}, Zipser calculates a Bondi mass loss formula by considering 
\begin{equation*}
\frac{\partial }{\partial u}m\left( t,u\right)
\end{equation*}%
with%
\begin{equation*}
\frac{\partial }{\partial u}m\left( t,u\right) =\frac{1}{2}\overline{%
atr\theta }m+\frac{r}{32\pi }\int_{S_{t,u}}a\left( \nabla _{N}\underline{\mu 
}+tr\theta \underline{\mu }\right) .
\end{equation*}%
As $l=a^{-1}\left( T+N\right) $ and $\underline{l}=a\left( T-N\right) $%
, 
\begin{eqnarray*}
a(\nabla _{N}\underline{\mu }+tr\theta \underline{\mu }) &=&\frac{1}{2}%
a^{2}\left( \mathbf{D}_{4}\underline{\mu }+tr\chi \underline{\mu }\right) \\
&&-\frac{1}{2}\left( \mathbf{D}_{3}\underline{\mu }+tr\underline{\chi }%
\underline{\mu }\right)
\end{eqnarray*}%
and 
\begin{eqnarray*}
\mathbf{D}_{4}\underline{\mu }+tr\chi \underline{\mu } &=&O(r^{-4}) \\
\mathbf{D}_{3}\underline{\mu }+tr\underline{\chi }\underline{\mu } &=&-\frac{%
1}{4}tr\chi \left \vert \widehat{\underline{\chi }}\right \vert ^{2}-\frac{1%
}{2}tr\chi \left \vert \underline{\alpha }\left( F\right) \right \vert
^{2}+O(r^{-4}).
\end{eqnarray*}%
Thus, 
\begin{equation*}
\frac{\partial }{\partial u}m\left( t,u\right) =\frac{r}{64\pi }%
\int_{S_{t,u}}tr\chi \left( \left \vert \widehat{\underline{\chi }}\right
\vert ^{2}+\frac{1}{2}\left \vert \underline{\alpha }\left( F\right) \right
\vert ^{2}\right) +O\left( r^{-1}\right) .
\end{equation*}%
Following \cite{sta}, Zipser uses the following facts to derive the Bondi mass loss 
formula: the metric $\widetilde{\gamma }=\phi _{t,u}^{\ast }\left(
r^{-2}\gamma \right) $ converges to the standard metric $\overset{\circ }{%
\gamma }$ of the unit sphere $S^{2}$ as $t\rightarrow \infty $ for each $u$ (%
$\phi _{t,u}^{\ast }$ is a diffeomorphism from $S^{2}$ to $S_{t,u}$), moreover $%
\frac{r}{2}tr\chi $ converges to $1$, and $r\widehat{\underline{\chi }}$
converges to $-2\Xi $. This yields  
\begin{equation*}
\frac{\partial }{\partial u}M\left( u\right) =\frac{1}{8\pi }%
\int_{S^{2}}\left( \left \vert \Xi \right \vert ^{2}+\frac{1}{2}\left \vert
A_{F}\right \vert ^{2}\right) d\mu _{\overset{\circ }{\gamma }}.
\end{equation*}%
The right-hand side of this expression is positive and integrable
in $u$. Thus, $M\left( u\right) $ is a non-decreasing function of $u$
and has finite limits $M\left( -\infty \right) $ for $u\rightarrow -\infty $
and $M\left( \infty \right) $ for $u\rightarrow \infty $. \ 
Further, Zipser concludes from (\ref{integrate mu}) 
that $M\left( -\infty \right) =0$, and $M\left( \infty \right) $ is the total mass.
\begin{The} \label{theBondimassformula1}
The Hawking mass $m\left( t,u\right) $ tends to the Bondi mass $M\left(
u\right) $ as $t \to \infty $ on any null hypersurface $C_{u}$. \
That is,  
\begin{equation*}
m(t,u)=M(u)+O(r^{-1}).
\end{equation*}
and $M\left( u\right) $ verifies the Bondi mass loss formula 
\begin{equation*}
\frac{\partial }{\partial u}M\left( u\right) =\frac{1}{8\pi }%
\int_{S^{2}}\left( \left| \Xi \right| ^{2}+\frac{1}{2}\left| A_{F}\right|
^{2}\right) d\mu _{\overset{\circ }{\gamma }}
\end{equation*}
with $d\mu _{\overset{\circ }{\gamma }}$ being the area element of the
standard unit sphere $S^{2}$.
\end{The}
We see that in the Bondi mass loss formula the limiting term $A_F$ of the
electromagnetic field comes in. 
At this point, let us compare this with the Bondi mass loss formula obtained in
\cite{sta} (p. 499): 
$\frac{\partial }{\partial u}M\left( u\right) =\frac{1}{8\pi }%
\int_{S^{2}} \left| \Xi \right| ^{2} d\mu _{\overset{\circ }{\gamma }}$. 
In fact, the electromagnetic field contributes to the change of the Bondi mass by 
$\frac{1}{16 \pi} \int_{S^{2}} \left| A_{F}\right|^{2} d\mu _{\overset{\circ
}{\gamma }}$.  \\ \\ 
The decay behavior of $A_F$ is the same as for $\Xi$. See theorems 
\ref{conclcurvandfieldcompts1} and \ref{conclXi1}. 
Similarly as in \cite{sta}, \cite{chrmemory} for the Einstein vacuum case, we can
define now the new function 
\be \label{en*1}
F = \frac{1}{8 } \int_{- \infty}^{+ \infty} \left( \mid \Xi \mid^2 + \frac{1}{2}
\mid A_F \mid^2  \right) du \ \ . 
\ee
Then $\frac{F}{4 \pi}$ is the total energy radiated to infinity in a given direction
per unit solid angle. 
Thus the integrand in (\ref{en*1}) is proportional to the power radiated to infinity
at a given retarded time $u$, in a given direction, 
per unit area on $S^2$ (per unit solid angle). Already in \cite{chrmemory}
Christodoulou tells us how to adapt the formula for $F$ when 
matter radiation is present, that is also in the EM case.  \\ \\
In the next two subsections, we also need the following theorem for $H$.
\begin{The}
The function $H$ satisfies
\begin{equation} \label{H1} 
\frac{\partial H}{\partial u} =0 
\end{equation}
\begin{equation}
\bar H = 0  \label{H2}
\end{equation}
\end{The}
{\bf Proof:} 
In the EV case, equation (\ref{H1}) is proved in conclusion 17.0.5 of \cite{sta} where one uses the fact that
\[ \nabla_N tr \chi' +\frac{1}{2}  \chi'  = O(r^{-3}). \]
In the EM case, it is easy to see that the additional terms involving the electromagnetic field are also $O(r^{-3})$.
Thus  equation (\ref{H1}) is still true in the EM case. \\ \\
In the EV case, equation (\ref{H2}) is proved in lemma 17.0.1 in \cite{sta}. In the proof, we need to show that 
$ r^2 \bar \delta  $ converges to $2 M(u)$. 
From Proposition 4.4.4 in \cite{sta}, we have
\[  4 \pi r^3 \bar \delta  = \int_{u_0}^u du'( \int_{S_{t,u}}  ar \hat \theta \cdot \hat \eta  - \frac{1}{2} \kappa(\delta - \bar \delta) - r a^{-1}  \nlap a \cdot \epsilon + r (div k)_N ) \]
Following the proof of lemma 17.0.1 in \cite{sta}, we see that 
\[   \int_{S_{t,u}}  ar \hat \theta \cdot \hat \eta  - \frac{1}{2} \kappa(\delta - \bar \delta) - r a^{-1}  \nlap a \cdot \epsilon  = r \int_{S^2} |\Xi|^2 d\mu _{\overset{\circ }{\gamma }} +O(1)\]
Moreover, in the EM case, due to the constraint equation, we have 
\begin{equation}  \label{dvikn} (div k)_N =  R_{0N}  = 2 F_{0}^{\ \rho} F_{N \rho} . \end{equation}
Using the equation (\ref{dvikn}), we see that
\[   \int_{S_{t,u}}  r (div k)_N = \frac{r}{2}\int_{S^2} |A_F|^2 d\mu _{\overset{\circ }{\gamma }} +O(1)\]
since the leading term of $F_{0A}$ and $F_{NA}$ are both $\frac{1}{2} \underline \alpha (F)_A$.
As a result, we can still conclude that 
\[ r \bar \delta  = \frac{2}{r^2} \int_{u_0}^u  r \frac{\partial }{\partial u}  m(t,u) +  O(r^{-1})\]
The rest of the proof follows easily from the proof of lemma 17.0.1 in \cite{sta}.
\subsection{Compare result with mass loss in Bondi coordinates}
In this subsection, we compare the mass loss formula obtained in \cite{zip},
\cite{zip2}, cited above in Theorem \ref{theBondimassformula1}, with the mass loss
formula in Bondi coordinates. 
Bondi coordinates are first defined by Bondi, van der Burg and Metzner in \cite{BBM}
for axially symmetric vacuum spacetimes. The main motivation for such coordinates is
to study gravitational radiation at null infinity. 
The form of the metric is chosen such that many computations are simplified at null
infinity. As a result, one can derive many useful theorems and formulae assuming the
existence of such coordinates. In particular, the Bondi mass loss formula in Bondi
coordinates is first derived for axially symmetric vacuum spacetimes in \cite{BBM}
and is 
generalized to the Einstein-Maxwell case in \cite{vanderBurg}.  However, for a given
spacetime, it is hard to tell whether such coordinates exist. On the other hand,
spacetimes studied in
\cite{zip} are obtained by evolving small initial data on a spacelike hypersurface
by EM equations. In particular, it is not clear whether all spacetimes studied in
\cite{zip} admit such coordinates. 
Here we show that for the leading term in the mass loss formula, the two different
coordinate systems give the same result. \\ \\
The mass loss formula from \cite{zip}, \cite{zip2} is:
\begin{equation*}
\frac{\partial }{\partial u}M\left( u\right) =\frac{1}{8\pi }%
\int_{S^{2}}\left( \left| \Xi \right| ^{2}+\frac{1}{2}\left| A_{F}\right|
^{2}\right) d\mu _{\overset{\circ }{\gamma }}.
\end{equation*}
Let us recall the mass loss formula and asymptotic expansion for solutions to the
Einstein-Maxwell equations in Bondi coordinates \cite{vanderBurg}. 
The line segment is:
$$
-UVdw^2-2Udwdr+\sigma_{ab}(dx^a+W^adw)(dx^b+W^bdw)\,\,\, a, b=2,3
$$
with the electromagnetic field given by a skew-symmetric two tensor $F_{\mu\nu}$. 
We have the following asymptotics for the line segment.
$$  V  =  1 -\frac{2m}{r} +O(r^{-2}) \text{, }  U  =  1 +O(r^{-2}) \text{ and } W^a
= O(r^{-2}).$$
\[ \sigma_{ab}=\left( \begin{array}{ccc}
r^2+2cr+ \cdots & -2dr \sin \theta +  \cdots \\
 -2dr \sin \theta +  \cdots & \sin^2 \theta(r^2-2cr) + \cdots
\end{array} \right).\] 
We also need the  following asymptotics for $F_{\mu\nu}$
\[ F_{w\theta}  = X + O(r^{-1})  \text{ and }F_{w\phi}  = Y\sin\theta  + O(r^{-1}).  \]
as well as
\[ F_{ra} = O(r^{-2}) \text{, } F_{ab} =O(1)  \text{ and } F_{wr} =O(r^{-2}).  \]
Level sets of $w$, $C_w$, are outgoing null hypersurfaces.
Each $C_w$ is then foliated by level sets of $r$, $S_{w,r}$. Hence, it is natural to
consider the following  pair of null vectors normal to $S_{w,r}$
$$e_3 = \frac{\partial}{\partial w} -W^c \frac{\partial}{\partial x^c} -\frac{V}{2}
\frac{\partial}{\partial r}  \text{ \, \, and \,  \,} e_4= \frac{\partial}{\partial
r} $$
since $e_4$ is a natural choice of null vector on $C_w$ and
\[  \lim_{r \to \infty} \langle e_3, e_4 \rangle= -1. \]
Let $M(w)$ be the Bondi mass. It is given by 
\[ M(w) = \frac{1}{8 \pi} \int_{S^2} m\, d\mu _{\overset{\circ }{\gamma }} .\] 
The mass loss formula reads 
$$
\frac{\partial }{\partial w}M(w) = - \int_{S^2} \left( (\partial_w c)^2 +(\partial_w
d)^2  +\frac{1}{2}(X^2+Y^2)  \right) d\mu _{\overset{\circ }{\gamma }}.
$$
\\ \\
To show that the two mass loss formulae agree, we prove that 
$$ |\Xi|^2 = (\partial_w c)^2 +(\partial_w d)^2 \text{ and } |A_F|^2 = X^2+Y^2. $$
First we compute
\begin{align*}
-\underline{\chi }^{\prime } (\frac{\partial}{\partial x^a} ,
\frac{\partial}{\partial x^b}) ={}& \langle   \frac{\partial}{\partial w} -W^c
\frac{\partial}{\partial x^c} -\frac{V}{2} \frac{\partial}{\partial r}  ,  \nabla
_{\frac{\partial}{\partial x^a}} \frac{\partial}{\partial x^b} \rangle  \\
={} & \langle \frac{\partial}{\partial w} -W^c \frac{\partial}{\partial x^c}
-\frac{V}{2} \frac{\partial}{\partial r}  , \Gamma_{ab}^r \frac{\partial}{\partial
r} +\Gamma_{ab}^w \frac{\partial}{\partial w} + \Gamma_{ab}^c
\frac{\partial}{\partial x^c}\rangle  \\ 
={} & \langle \frac{\partial}{\partial w} -W^c \frac{\partial}{\partial x^c}
-\frac{V}{2} \frac{\partial}{\partial r}  , \Gamma_{ab}^r \frac{\partial}{\partial
r} +\Gamma_{ab}^w \frac{\partial}{\partial w} \rangle  \\ 
\end{align*}
The Christoffel symbols are 
\begin{align*} \Gamma^w_{ab}  ={} &   -\frac{1}{2}     g^{wr}  \partial_r g_{ab} =
\frac{1}{2} \partial_r g_{ab}  +O(1)\\
\Gamma^r_{ab}  ={} &    \frac{1}{2}     g^{wr} (  \partial_b g_{wa}+ \partial_a
g_{wb} -  \partial_w g_{ab}) +
\frac{1}{2}     g^{rr} ( -  \partial_r g_{ab}) + \frac{1}{2}     g^{rc} ( 
\partial_b g_{ca}+ \partial_a g_{cb}-  \partial_c g_{ab} )\\
={}& \frac{1}{2}\partial_w g_{ab}-\frac{1}{2} \partial_r g_{ab}  +O(1)
 \end{align*}
As a result, one finds 
$$ \langle   \frac{\partial}{\partial w} -W^c \frac{\partial}{\partial x^c}
-\frac{V}{2} \frac{\partial}{\partial r}  ,  \nabla _{\frac{\partial}{\partial x^a}}
\frac{\partial}{\partial x^b} \rangle  
=  \frac{1}{4}\partial_r \sigma_{ab} - \frac{1}{2} \partial_w \sigma_{ab} +O(1).$$
One can easily see that up to $O(1)$ terms,  $\partial_w \sigma_{ab}$ is traceless 
and  $\partial_r \sigma_{ab}$ has zero traceless part. As a result, 
$$  | \Xi|^2 = (\partial_w c)^2 +(\partial_w d)^2. $$
For the second equality, we use the expression for $e_3$ and the asymptotics for
$F_{\mu \nu}$. A direct computation shows that 
$$|A_F|^2 = X^2+Y^2. $$ 
\subsection{Permanent Displacement Formula}
\label{permanentdispl***}
The permanent displacement of the test masses of a laser interferometer
gravitational-wave detector is 
governed by $\Sigma^+ - \Sigma^-$. Christodoulou showed in \cite{chrmemory} how this
works. 
We discuss the corresponding wave experiment in the EM case in section \ref{wave}. 
In the following, we are going to state and prove a theorem for $\Sigma^+ -
\Sigma^-$ in the EM case. 
We point out that the final formula - even though in its form identical to the one
obtained by Christodoulou and Klainerman in \cite{sta} - 
differs from the EV case by a contribution from the electromagnetic field. 
The form of the formula is not altered due to 
the fact that the corresponding extra electromagnetic terms cancel. However, the
limiting term $A_F$ enters the 
new formula nonlinearly. \\ 
\begin{The} \label{displ*1}
Let 
$\Sigma^+ (\cdot) = \lim_{u \to \infty} \Sigma (u, \cdot)$ and 
$\Sigma^- (\cdot) = \lim_{u \to - \infty} \Sigma (u, \cdot)$. 
Let 
\be \label{Thm*FXiAF*1}
F (\cdot)  =   \int_{- \infty}^{\infty} 
\big( 
\mid \Xi (u, \cdot) \mid^2 + \frac{1}{2} \mid A_F (u, \cdot) \mid^2  
\big)
du  \ \ . 
\ee
Moreover, let 
$\Phi$ be the solution with $\bar{\Phi} = 0$ on $S^2$ of the equation 
\[
\stackrel{\circ}{\slap} \Phi = F - \bar{F}   \ \ . 
\]
Then 
$\Sigma^+ - \Sigma^-$ is given by the following equation on $S^2$: 
\be \label{Thm*divSigma+-*2}
\stackrel{\circ}{\dlap} (\Sigma^+ - \Sigma^-) = \stackrel{\circ}{\nlap} \Phi \ \ . 
\ee
\end{The}
{\bf Proof:}  
Equation (\ref{Sigmau*1}) in theorem \ref{conclXi1} yields that 
$\Sigma$ tends to limits $\Sigma^+$ as $u \to \infty$ and $\Sigma^-$ as $u \to -
\infty$. Also, it is 
\[
\Sigma (u) = \Sigma^- - \int_{- \infty}^u \Xi (u') du' 
\]
and 
\[
\Sigma^+ - \Sigma^- = - \int_{- \infty}^{\infty} \Xi (u') du'  \ \ . 
\]
When taking the limits for the 
Hodge system ((\ref{systeps1}), (\ref{systeps2})) on $C_u$ as $t \to \infty$, 
we will compute the corresponding limits for and involving $\Psi$, $\Psi'$. 
From Zipser's work \cite{zip2}, chapter 9.1, (9.13) and lemma 9, we know 
\bea
\triangle \Psi & = & r \mid \hat{\eta} \mid^2  - \frac{r}{4} \mid \underline{\alpha}
(F) \mid^2  \label{trianglePsi1}  \\ 
\triangle \Psi' & = & - r a^{-1} \lambda \big( \mid \hat{\eta} \mid^2 -
\overline{\mid \hat{\eta} \mid^2} \big) 
 + \frac{r^2 a^{-1}}{4} \big( a  \Dlap_4 \mid \underline{\alpha} (F) \mid^2  -
\overline{a  \Dlap_4 \mid \underline{\alpha} (F) \mid^2}  \big) 
\label{trianglePsip1}
\eea
whereas in the work of Christodoulou and Klainerman \cite{sta}, chapter 11.2,
(11.2.2b) and (11.2.7b)
it is 
\bea
\triangle \Psi & = & r \mid \hat{\eta} \mid^2  \label{trianglePsi2} \\ 
\triangle \Psi' & = & - r a^{-1} \lambda \big( \mid \hat{\eta} \mid^2 -
\overline{\mid \hat{\eta} \mid^2} \big)  \label{trianglePsip2}  \ \ . 
\eea
We compute the following limits in the new case. 
\bea 
\lim_{C_u, t \to \infty} \Psi = \mathbf{\Psi} \ \ \ \ \quad & & \ \ \ \ \quad 
\lim_{C_u, t \to \infty} \Psi' = \mathbf{\Psi'}  \nonumber \\ 
\lim_{C_u, t \to \infty} r \nabla_N \Psi = \Omega(u, \cdot )  \ \ & & \ \ 
\lim_{C_u, t \to \infty} r \nabla_N \Psi' = \Omega'(u, \cdot ) 
\label{limitsPsiOmega*1}.
\eea 
We proceed by investigating the Hodge system for $\epsilon$. 
The Hodge system for $\epsilon$ reads, see also \cite{zip2}, chapter 9: 
\bea
\dlap \epsilon \ & = & \ 
- \nabla_N \delta - \frac{3}{2} tr \theta \delta + \hat {\eta} \cdot \hat{\theta} 
\nonumber \\ 
\ & & \ 
- 2 (a^{-1} \nlap a ) \cdot \epsilon + \frac{1}{4} \mid \alpha (F) \mid^2 -
\frac{1}{4} 
\mid \underline{\alpha} (F) \mid^2 \label{systeps1} \\ 
\clap \   \ \epsilon  \ & = & \ 
\sigma (W) + \hat{\theta} \wedge \hat{\eta} \ \ .  \label{systeps2}
\eea
We observe that the $\clap \ $ equation coincides with the one obtained by
Christodoulou and Klainerman in \cite{sta}, whereas 
the $\dlap $ equation contains the extra terms $\mid \alpha (F) \mid^2$ and $\mid
\underline{\alpha} (F) \mid^2$ from 
the electromagnetic field. (See \cite{sta}, chapter 17, ((17.0.12a), (17.0.12b)).) 
According to \cite{zip2}, chapter 9, one has: 
\bea
\nabla_N \delta - \hat{\theta} \cdot \hat{\eta} + \frac{1}{4} \mid
\underline{\alpha} (F) \mid^2 \ & = & \ 
- 2r^{-3} (\nabla_N r) p + r^{-2} \nabla_N p - r^{-2} (\nabla_N r) \nabla_N \Psi +
r^{-1} \nabla_N^2 \Psi  \nonumber \\ 
\ & = & \ 
- \hat{\chi} \cdot \hat{\eta} - r^{-1} \slap \Psi - r^{-2} 
\big( 
r tr \theta + a^{-1} \lambda \big) \nabla_N \Psi   \nonumber  \\ 
\ & & \ 
- r^{-1} a^{-1} \nlap a \cdot \nlap \Psi + r^{-2} \nabla_N p - 2 r^{-3} a^{-1}
\lambda p  \label{nablaNdelta*1}
\eea
with 
\[
p = r \nabla_N q + q' + \Psi' 
\]
This differs from \cite{sta}, chapter 17, (17.0.12c) by the extra curvature term from
the electromagnetic field. 
Zipser derived in \cite{zip2}, chapter 9, Lemma 9, 
\be \label{triangleq*1}
\triangle q \ = \ r (\mu - \overline{\mu}) +I 
\ee
where
\beas
I \ & = & \ \frac{1}{2} \ ^{(rN)} \hat{\pi}_{ij} k_{ij} + \frac{r}{4} \mid \alpha
(F) \mid^2 - \frac{r}{4} \mid \underline{\alpha} (F) \mid^2 
- \triangle \Psi \\ 
\ & = & \ 
r \hat{\chi} \cdot \hat{\eta} - \kappa \delta - 2 r a^{-1} \nlap a \cdot \epsilon +
\frac{r}{4} \mid \alpha (F) \mid^2 \  
\eeas
and $\mu$ is the mass aspect function given by 
\[  \mu=  - \rho(W) - \widehat \chi   \cdot \hat \eta \ .\]
Recall the radial decomposition of $\triangle$ to be 
$\nabla_N^2 = \triangle -  tr \theta \nabla_N - \slap - a^{-1} \nlap a \cdot \nlap$. 
Now, we obtain from the last equations that 
\bea 
\triangle q \ & = & \ 
\nabla_N^2 q +  tr \theta \nabla_N q + \slap q + a^{-1} \nlap a \cdot \nlap q  
\nonumber \\ 
\ & = & \ 
-r(\rho - \bar{\rho}) - r \overline{\hat{\chi} \cdot \hat{\eta}} - \kappa \delta - 2
r a^{-1} \nlap a \cdot \epsilon 
+ \frac{r}{4} \mid \alpha (F) \mid^2  \label{triangle*2} 
\eea
We proceed as follows: 
Substituting first for $\nabla_N p$ from (\ref{triangle*2}) in (\ref{nablaNdelta*1})
and then 
the resulting terms from (\ref{nablaNdelta*1}) in (\ref{systeps1}) yields 
\bea 
\dlap \epsilon \ & = & \ 
\rho - \bar{\rho} + \hat{\chi} \cdot \hat{\eta} - \overline{\hat{\chi} \cdot
\hat{\eta}} 
+ r^{-1} \slap \Psi - r^{-2} \nabla_N \Psi' - r^{-3} a^{-1} \lambda \Psi' + l.o.t. 
\label{systeps3} \\ 
\clap \   \ \epsilon  \ & = & \ 
\sigma (W) + \hat{\theta} \wedge \hat{\eta}  \label{systeps4}
\eea
Let 
\[ E= \lim_{C_{u},t\rightarrow \infty }\left( r \epsilon \right ) \]
We multiply equations (\ref{systeps3}) and  (\ref{systeps4}) by $r^3$ and take the limits on $C_u$ as $t \to \infty$. This yields: 
\bea 
\stackrel{\circ}{\clap} \ \ E \ & = & \ Q + \Sigma \wedge \Xi  \label{systeps*5*1}  \\ 
\stackrel{\circ}{\dlap} E \ & = & \ P - \bar{P} + \Sigma \cdot \Xi - 
\overline{\Sigma \cdot \Xi}  \nonumber \\ 
\ & & \ 
+ \stackrel{\circ}{\slap} \Psi - \Psi' - \Omega' \ \ .  \label{systeps*5}
\eea 
Then we investigate the limits as $u \to + \infty$ and $u \to - \infty$. 
Considering the last equations for $\epsilon$, respectively $E$, and using theorems 
\ref{conclXi1} and 
\ref{conclcurvandfieldcompts1} 
one finds that $E$ tends to a limit $E^+$ as $u \to + \infty$ and to $E^-$ as $u \to
- \infty$. \\ \\ 
By conclusions along the lines of \cite{sta}, chapter 17, 
we obtain 
\[
\stackrel{\circ}{\clap} \ \ (E^+ - E^-  ) \  =  \ 0 \\ 
\]
In order to compute $\stackrel{\circ}{\dlap} (E^+ - E^-  )$, we have to consider
especially the corresponding limits for the terms 
involving $\mathbf{\Psi}$ and $\mathbf{\Psi'}$, that is also $\Omega'$. \\ \\ 
Much like Christodoulou and Klainerman computed the formulas in lemma 17.0.2, on
page 504 of \cite{sta}, we derive the new results 
in which the electromagnetic field term $\underline{\alpha}(F)$, respectively its
limit $A_F$, is present. 
To do that, we use the fact that 
\be
\Dlap_4 \underline{\alpha} (F)_A = - \frac{1}{2} tr \chi \underline{\alpha} (F)_A +
l.o.t.  \label{Dlap4alphaundF*1}
\ee
The derivation of (\ref{Dlap4alphaundF*1}) can be found in Zipser's work \cite{zip2}
on p. 351 formula (4.15). 
Now, considering (\ref{trianglePsip1}) and using (\ref{Dlap4alphaundF*1}), 
we find that 
\be
\Dlap_4 \mid \underline{\alpha} (F) \mid^2 = -  tr \chi \mid \underline{\alpha} (F)
\mid^2 + l.o.t.  \label{Dlap4alphaundF**}
\ee
Using (\ref{Dlap4alphaundF**}), (\ref{trianglePsi2}) and (\ref{trianglePsip2}), we
deduce formulas for 
$\mathbf{\Psi}$, $\mathbf{\Psi'}$, $\Omega$, $\Omega'$ by computing the limits
(\ref{limitsPsiOmega*1}). 
We give the formulas:  
\beas
\mathbf{\Psi} & = & 
- \frac{1}{2^{\frac{1}{2}} 4 \pi} \int_{- \infty}^{+ \infty} \Big{ \{ } 
\int_{S^2} \frac{\mid \Xi \mid^2 (u', \omega' ) }{( 1 - \omega \omega'
)^{\frac{1}{2}}} 
d \omega' 
+
\frac{1}{2} 
\int_{S^2} \frac{\mid A_F \mid^2 (u', \omega' ) }{( 1 - \omega \omega'
)^{\frac{1}{2}}} 
d \omega' 
\Big{ \} } 
d u' \\ 
\mathbf{\Psi'} & = & 
\frac{1}{2^{\frac{1}{2}} 4 \pi} \int_{- \infty}^{+ \infty} \Big{ \{ } 
\int_{S^2} \frac{\mid \Xi \mid^2 (u', \omega' ) - \overline{\mid \Xi \mid^2 (u')} 
}{( 1 - \omega \omega' )^{\frac{1}{2}}} 
d \omega' 
+ 
\frac{1}{2} 
\int_{S^2} \frac{\mid A_F \mid^2 (u', \omega' )  - \overline{\mid A_F \mid^2 (u')} 
}{( 1 - \omega \omega' )^{\frac{1}{2}}} 
d \omega' 
\Big{ \} } 
d u' \\ 
\Omega & = & 
\frac{1}{2^{\frac{3}{2}} 4 \pi} \int_{- \infty}^{+ \infty} \Big{ \{ } 
\int_{S^2} \frac{\mid \Xi \mid^2 (u', \omega' ) }{( 1 - \omega \omega'
)^{\frac{1}{2}}} 
d \omega' 
+ 
\frac{1}{2} 
\int_{S^2} \frac{\mid A_F \mid^2 (u', \omega' ) }{( 1 - \omega \omega'
)^{\frac{1}{2}}} 
d \omega' 
\Big{ \} } 
d u'     \\ 
& & \ \ 
+ \frac{1}{2} 
 \int_{- \infty}^{+ \infty} 
 \Big{ \{ } 
sgn(u - u') 
\Big( 
\mid \Xi \mid^2 (u', \omega') + \frac{1}{2} \mid A_F \mid^2 (u', \omega') 
\Big)
\Big{ \} } du' \\ 
\Omega' & = & 
- \frac{1}{2^{\frac{3}{2}} 4 \pi} \int_{- \infty}^{+ \infty} \Big{ \{ } 
\int_{S^2} \frac{\mid \Xi \mid^2 (u', \omega' ) - \overline{\mid \Xi \mid^2 (u')}
}{( 1 - \omega \omega' )^{\frac{1}{2}}} 
d \omega' 
+ 
\frac{1}{2} 
\int_{S^2} \frac{\mid A_F \mid^2 (u', \omega' ) - \overline{\mid A_F \mid^2 (u')}
}{( 1 - \omega \omega' )^{\frac{1}{2}}} 
d \omega' 
\Big{ \} } 
d u'     \\ 
& & \ \ 
- \frac{1}{2} 
 \int_{- \infty}^{+ \infty} 
 \Big{ \{ } 
sgn(u - u') 
\Big( \Big(
\mid \Xi \mid^2 (u', \omega') - \overline{\mid \Xi \mid^2 (u')} \Big)
+ \frac{1}{2} \Big( \mid A_F \mid^2 (u', \omega') -  \overline{\mid A_F \mid^2 (u')} 
\Big)
\Big)
\Big{ \} } du' \\ 
\eeas
Straightforward calculation shows that when evaluating the difference of the limits
as \\ 
$u \to + \infty$ and $u \to - \infty$ 
in (\ref{systeps*5}), 
the contribution of 
$\stackrel{\circ}{\slap} \Psi$, $\Psi'$ and $\Omega'$ comes only from terms in
$\Omega'$. 
We find that $\Omega'$ tends to limits $\Omega'^+ (\cdot)$ and $\Omega'^- (\cdot)$
as $t \to \infty$ and $t \to - \infty$, respectively. 
Thus, we conclude 
\be
\Omega'^+ (\cdot) - \Omega'^- (\cdot) \ = \ 
\int_{- \infty}^{+ \infty} \big( \ 
\mid \Xi (u, \cdot) \mid^2 - \overline{\mid \Xi (u, \cdot) \mid^2} 
+ \frac{1}{2} \mid A_F (u, \cdot) \mid^2 - \frac{1}{2} \overline{\mid A_F (u, \cdot)
\mid^2} 
\  \big) du \ \ . 
\ee
Finally, we obtain 
\bea
\stackrel{\circ}{\dlap} (E^+ - E^-  ) \ & = & \ - \Omega'^+  + \Omega'^-   \\ 
\ & = & \ 
\int_{- \infty}^{+ \infty} \big( \ -
\mid \Xi (u, \cdot) \mid^2 + \overline{\mid \Xi (u, \cdot) \mid^2} 
- \frac{1}{2} \mid A_F (u, \cdot) \mid^2 + \frac{1}{2} \overline{\mid A_F (u, \cdot)
\mid^2} 
\  \big) du \ \ .  \nonumber 
\eea
Proceeding along the lines of \cite{sta}, chapter 17, it is 
\be
(E^+ - E^- ) \ = \ \stackrel{\circ}{\nlap} \Phi 
\ee
with $\Phi$ being the solution of vanishing mean of 
\[
\stackrel{\circ}{\slap} \Phi \ = \ - \Omega'^+  + \Omega'^-   \ \ \mbox{ on } S^2 \
\ . 
\]
Accordingly, also by conclusions along the lines of \cite{sta}, chapter 17, we
derive (\ref{*circSigma**1}). 
To see this, we consider the normalized null Codazzi equation 
\be \label{*normalizednullcodazzi**1}
(\dlap \hat{\chi})_A  - \frac{1}{2} \nlap_A tr \chi + \epsilon_B \hat{\chi}_{AB} -
\frac{1}{2} \epsilon_A tr \chi  =  
- \beta (W)_A - \rho(F) \alpha(F)_A - \epsilon_{AB} \sigma(F) \alpha(F)_B  
\ee
Multiply equation (\ref{*normalizednullcodazzi**1}) by $r^3$ and take the limit as 
$t \to \infty$ on $C_u$. We obtain  
\[ \stackrel{\circ}{\dlap} \Sigma = \stackrel{\circ}{\nlap} H+E \]
as in \cite{sta} p. 510, conclusion 17.0.8 since the extra terms from the electromagnetic field in (\ref{*normalizednullcodazzi**1})
decay fast enough.
Due to equation (\ref{H1}), we conclude 
\be \label{*circSigma**1}
\stackrel{\circ}{\dlap} (\Sigma^+ - \Sigma^-) = E^+ - E^- \ \ . 
\ee
Thus, the theorem is proven.
\subsection{Limit for $r$ as $t \to \infty$ on Null Hypersurface $C_u$}
We shall use the fact that the constraint on the spacelike scalar
curvature, which is given by 
\begin{equation*}
R=\left \vert k\right \vert ^{2}+R_{00},
\end{equation*}%
differs from the constraint in the vacuum case only by the term $R_{00}$,
which is a quadratic in $F$.  \\ \\ 
Building on the results of D. Christodoulou and S. Klainerman in \cite{sta} as well
as the results of N. Zipser in \cite{zip2} (i.e. \cite{zip}), 
we can now prove the following results. 
\begin{The} \label{r*1}
As $t \to \infty$ we obtain on any null hypersurface $C_u$ 
\[
r = t - 2 M(\infty) \log t + O(1)  \ \ . 
\]
\end{The}
{\bf Proof:} 
We recall from \cite{sta}, p. 503, with $\phi' = \phi -1$, 
\beas
\frac{dr}{dt} & = & \frac{r}{2} \overline{\phi tr \chi'} \\ 
& = &  
 \frac{r}{2} \overline{(1 + \phi')(\frac{2}{r} + (tr \chi' - \frac{2}{r}))} \\ 
& = & 
1 + \overline{\phi'} + O(r^{-2})  
\eeas 
In the last equality, we use equation (\ref{H2}).

From \cite{zip2}, p. 465 
\be \label{R00*1}
R_{00} = \frac{1}{2} ( \mid \underline{\alpha}(F) \mid^2 + \mid \alpha(F) \mid^2  )
+ \rho(F)^2 + \sigma(F)^2 
\ee
with $\underline{\alpha}(F),  \alpha(F),  \rho(F), \sigma(F)$ the components of the
electromagnetic field.  Moreover, 
the lapse equation in our situation is given by 
\be \label{lapset*1}
\triangle \phi = ( \mid k \mid^2 + R_{00} ) \phi \ \ . 
\ee
We integrate the lapse equation (\ref{lapset*1}) on $H_t$ in the interior of $S_{t,
u'}$ to obtain 
\[
\int_{S_{t, u}} \nabla_N \phi' = \int_{u_0}^u du' \int_{S_{t, u'}} a \phi ( \mid k
\mid^2  + R_{00} ) \ \ . 
\]
In view of (\ref{R00*1}) and the fact that all the terms on the right hand side of
(\ref{R00*1}) except $\underline{\alpha}(F)$ are of 
lower order, we estimate 
\beas
\int_{S_{t, u}} \nabla_N \phi'  & = &   \int_{u_0}^u du' \int_{S_{t, u'}} a \phi (
\mid k \mid^2  + \frac{1}{2}  \mid \underline{\alpha}(F) \mid^2  )  + l.o.t.  \\ 
\eeas
We see that 
\[
\int_{S_{t, u'}}  a \phi ( \mid k \mid^2  + \frac{1}{2}  \mid \underline{\alpha}(F)
\mid^2  )  \to 
\int_{S^2} \mid \Xi \mid^2 + \frac{1}{2} \mid A_F \mid^2 \ \ . 
\]
Consider the Bondi mass loss formula in theorem \ref{theBondimassformula1}. Then, as
$t \to \infty$ we conclude 
\be \label{nablaNM**1}
\int_{S_{t,u}} \nabla_N \phi' - 8 \pi M(u) = O(r^{-1}) 
\ee
on each $C_u$. 
In view of $\phi'$ we compute: 
\beas
\overline{\phi'} & = & \frac{1}{4 \pi r^2} \int_{S_{t,u}} \phi' = - \frac{1}{4 \pi}
\int_{B} div (r^{-2}  \phi' N) \\ 
& = & 
\frac{1}{4 \pi}  \int_{B} \big( 
- \frac{1}{a(r(t,u'))^2} \overline{a tr \theta} N \phi' 
+ \frac{1}{(r(t,u'))^2} \phi' \underbrace{(div N)}_{= tr \theta} 
+ \frac{1}{(r(t,u'))^2} \nabla_N \phi' 
 \big) \\ 
 & = & 
- \frac{1}{4 \pi} \int_u^{\infty} \frac{1}{(r(t,u'))^2}  du' 
\big(
 \int_{S_{t, u'}} 
 a \nabla_N \phi' + (a tr \theta - \overline{a tr \theta}) \phi'  
  \big) \\ 
 & = & 
 - \frac{1}{4 \pi} \int_u^{\infty} \frac{1}{(r(t,u'))^2}  du' 
 \big( 
  \int_{S_{t, u'}}  \nabla_N \phi'  
 \big) 
 + O(r^{-2}) \ \ . 
\eeas
where $B$ denotes the exterior of $S_{t,u}$.
Therefore, from (\ref{nablaNM**1}) it follows on $C_u$ as $t \to \infty$, 
\[
\overline{\phi'} (t, u) = - 2 \int_u^{\infty}  \frac{1}{(r(t,u'))^2} M(u') du' +
O(r^{-2}) 
= - \frac{2}{r} M(\infty) + O(r^{-2}) \ \ . 
\] 
Thus, we obtain on any cone $C_u$ for $t \to \infty$, 
\be
\frac{dr}{dt} = 1 - \frac{2}{r} M(\infty) + O(r^{-2}) \ \ . 
\ee
Thus, the statement of our theorem follows, which closes the proof. \\ \\ 
%
%
%
%
%
%
%
\section{Wave Experiments}
\label{wave}
We are now going to show how the results above relate to experiment. 
In \cite{chrmemory} Christodoulou established his breaking result on the nonlinear
memory effect. 
The idea of the gravitational wave experiment and setup is given in
\cite{chrmemory}, discussing 
a laser interferometer gravitational-wave detector. 
There Christodoulou explained how the theoretical result on $\Sigma^+ - \Sigma^-$ 
leads to an effect measurable by such detectors. 
This effect manifests itself in a permanent displacement of the 
test masses 
of the detector after a wave train has passed. 
In the present EM case, we find a result on the displacement of test masses which is
twofold. 
Considering the Jacobi equation (see (\ref{Jac***3})), the highest order term
remains unchanged. 
However, there is an extra term at highest order from the electromagnetic field in
the formula for 
$\Sigma^+ - \Sigma^-$, the permanent displacement of test masses, as we have shown
in the proof of theorem \ref{displ*1}. 
In the present chapter, we shall show how the  electromagnetic field enters the
experiment and we will derive results for this case.  \\ \\ 
We will follow the lines of argumentation by Christodoulou in \cite{chrmemory} and
\cite{chrdlbmathpgrt}. \\ \\  
Let us briefly review the setup of a laser interferometer experiment: 
Three test masses are suspended by pendulums of equal length.  
Denote by $m_0$ the reference mass, which is also the location of the beam splitter. 
For timelike scales 
much shorter than the period of the pendulums, the motion of the masses in the
horizontal plane can be considered free. 
By laser interferometry the distance of the masses $m_1$ and $m_2$ from the
reference mass $m_0$ is measured. 
Whenever the light travel times between $m_0$ and $m_1$ and $m_2$, respectively, 
differ, this shows in a difference of phase of the laser light at $m_0$. \\ \\ 
The masses $m_0$, $m_1$, $m_2$ move along geodesics 
$\gamma_0$, $\gamma_1$, $\gamma_2$ in spacetime.  
Let $T$ be the unit future-directed tangent vectorfield of $\gamma_0$ and $t$ 
the arc length along $\gamma_0$. 
For each $t$ denote by $H_t$ the spacelike, geodesic hyperplane through 
$\gamma_0(t)$ orthogonal to $T$.  \\ \\ 
Take $(E_1, E_2, E_3)$ to be an orthonormal frame for $H_0$ 
at $\gamma_0 (0)$, and parallelly propagate it along $\gamma_0$ to obtain 
the orthonormal frame field 
$(T, E_1, E_2, E_3)$ along $\gamma_0$. 
It follows that 
$(E_1, E_2, E_3)$ at each $t$ is an orthonormal frame for 
$H_t$ at $\gamma_0 (t)$. 
Then one assigns to a point $p$ in spacetime, 
lying in a neighbourhood of 
$\gamma_0$, the 
cylindrical normal coordinates 
$(t, x^1, x^2, x^3)$, based on $\gamma_0$, if 
$p \in H_t$ and $p = exp X$ with 
$X = \sum_i x^i E_i \in T_{\gamma_0(t)} H_t$. 
Denote by $d$ the distance of $p$ from the center $\gamma_0(t)$ on $H_t$, that is, 
$d  =  \mid X \mid  =  \sqrt{\sum_i  (x^i)^2}$. 
The difference between the metric and the Minkowski metric $\eta_{\mu \nu}$ in these
coordinates is: 
\be
g_{\mu \nu} \ - \ \eta_{\mu \nu} \ = \ O \ (R  \ d^2) \ \ . 
\ee
As usual, we put $c=1$. Now, let $\tau$ be 
the time scale in which the curvature 
varies significantly. 
Then, the 
displacements of the masses from their 
initial positions will be  
$O (R \tau^2)$. 
Assume that  
\be \label{d/tau<1}
\frac{d}{\tau } \ <  < \ 1 \ \ . 
\ee
Then we can read off from the differences in phase 
of the laser light 
differences in distance of $m_1$ and $m_2$ from $m_0$. 
Further, 
in view of (\ref{d/tau<1}) one can 
replace the geodesic equation for $\gamma_1$ and $\gamma_2$ by the Jacobi 
equation (geodesic deviation from $\gamma_0$).
\be \label{Jacobi**1}
\frac{d^2 x^k}{d t^2} \ = \ - \ R_{kTlT} \ x^l 
\ee
with 
$R_{kTlT} \ = \ R \ (E_k, T, E_l, T)$. 
One is free to assume that the source is in the $E_3$-direction. 
This was derived by Christodoulou for the EV case in \cite{chrmemory}, and can also
be found in his 
\cite{chrdlbmathpgrt}. \\ \\ 
We now investigate the formula (\ref{Jacobi**1}) for the Einstein-Maxwell situation.
If we assume the test masses not to be 
charged, then formula (\ref{Jacobi**1}) stays the same, but through the EM equations
and in view of 
(\ref{Riem*tracelessRicci}) the electromagnetic field comes in. We shall see that it
enters at lower order though. 
From (\ref{Riem*tracelessRicci}) one can write 
\bea
R_{k0l0} &=&W_{k0l0}+\frac{1}{2}(g_{kl}R_{00}+
g_{00}R_{kl}-g_{0l }R_{k0 }-g_{k0 }R_{0l})    \ \ .  \label{RiemtracelessRicciEM*1}
\eea
As there is from the EM equations: 
\[
R_{00} = 8 \pi T_{00} \ \ , 
\]
and in particular, we have 
\be
R_{00} = \frac{1}{2} ( \mid \underline{\alpha} (F) \mid^2 + \mid \alpha (F) \mid^2 
) + \rho(F)^2 + \sigma(F)^2 
\ee
we can investigate the components of the Ricci curvature on the right hand side of
(\ref{RiemtracelessRicciEM*1}). 
The component $R_{00}$ includes the term $\mid \underline{\alpha}(F) \mid^2$. Recall
that 
$\underline{\alpha}(F)$ is the part of the 
electromagnetic field with worst decay behavior. However it enters as a quadratic
the formula for $R_{00}$. \\ \\ 
To proceed, we consider 
$ L = T - E_3$, $\underline{L} = T + E_3$. 
The leading components of the curvature are  
\bea
\underline{\alpha}_{AB} (W) \ & = & \ R \ (E_A, \ \underline{L}, \ E_B, \
\underline{L}) \\ 
\underline{\alpha}_{AB} (W) \ & = & \ \frac{A_{AB}(W)}{r} \ + \ o \ (r^{-2}) \ . 
\eea
And the leading components of the electromagnetic field are 
\bea
\underline{\alpha}_{A} (F) \ & = & \ F (E_A, \underline{L}) \  \\ 
\underline{\alpha}_{A} (F) \ & = & \ \frac{A_{A}(F)}{r} \ + \ o \ (r^{-2}) \ . 
\eea
In the following the $k$th Cartesian coordinate of the mass $m_A$ for $A = 1,2$ will
be denoted by 
$x^k_{\ (A)}$. 
Then the Jacobi equation 
becomes 
\be \label{Jac***3}
\frac{d^2 \ x^k_{\ (A)}}{d \ t^2} \ =  \ - \ \frac{1}{4} \ r^{-1}  \ A_{AB} \ x^l_{\
(B)} 
\ - \ \frac{1}{8} r^{-2} \mid A_F \mid^2 x^l_{\ (B)} 
\ + \ 
O \ (r^{-2}) 
\ee
that is 
\bea
\frac{d^2 \ x^3_{\ (C)}}{d \  t^2} \ & =  & \ 0 \\ 
\frac{d^2 \ x^A_{\ (C)}}{d \ t^2} \ & =  & \  - \ \frac{1}{4} \ r^{-1}  \ A_{AB} \
x^B_{\ (D)} 
\ - \ \frac{1}{8} r^{-2} \mid A_F \mid^2 x^B_{\ (D)} 
\ + \ 
O \ (r^{-2}) \ . 
\eea
From the Jacobi equation (\ref{Jac***3}) we see that the electromagnetic field 
enters on the right hand side at order $(r^{-2})$ only. 
Thus, we have shown that the electromagnetic field does not contribute at leading
order to the deviation measured 
by the Jacobi equation. 
Therefore, at leading order, we can rely on the results for the Einstein vacuum
case, derived by Christodoulou in 
\cite{chrmemory}. Instead of (\ref{Jac***3}) he obtained 
\be \label{Jac***4}
\frac{d^2 \ x^k_{\ (A)}}{d \ t^2} \ =  \ - \ \frac{1}{4} \ r^{-1}  \ A_{AB} \ x^l_{\
(B)} \ + \ 
O \ (r^{-2}) 
\ee
As in \cite{chrmemory} one obtains 
that in the vertical direction there is no acceleration to leading order $(r^{-1})$. 
Initially $m_1$ and $m_2$ are at rest at equal distance $d_0$ and at right angles from 
$m_0$. 
This implies the following initial conditions, as $t \to - \infty$:
$x^3_{\ (A)}  =  0 \ , \ \dot{x}^3_{\ (A)}  = 0 \ , \ 
x^B_{\ (A)} =  d_0  \delta^B_A  \ , \ 
\dot{x}^B_{\ (A)}  =  0$. 
The right hand side being very small, one can substitute the initial values on the
right hand side. 
Then the motion is confined to the horizontal plane. One has to leading order: 
\be
\stackrel{\cdot \cdot}{x}^A_{\ (B)} \ = \ - \ \frac{1}{4} \ r^{-1} \ d_0 \ A_{AB}  \
\ . 
\ee
One obtains 
\be 
\dot{x}^A_{\ (B)} \ (t) \ = \ - \ \frac{1}{4} \  d_0 \  r^{-1} \ 
\int_{- \infty}^t \ A_{AB} \ (u) \ d u \ . 
\ee
In the following, let us revisit the result (\ref{reschrmem****4}) from
\cite{chrmemory}. 
In view of equation (\ref{Xiu*1}), i.e. 
$
\frac{\partial \Xi }{\partial u}  = -  \frac{1}{4} \ A $
and 
$\lim_{\mid u \mid \to \infty}  \Xi  =  0$ 
we obtain 
\be
- \ \int_{- \infty}^t \ A_{AB} \ (u) \ d u \ = \ \Xi \ (t) \ 
\ee
and 
\be
\dot{x}^A_{\ (B)} \ (t) \ = \ \frac{d_0}{r}  \ \Xi_{AB} \ (t) \ . 
\ee
As $\Xi \to 0$ for $u \to \infty$, the test masses return to rest after 
the passage of the gravitational wave. 
Taking into account (\ref{Sigmau*1}), i.e. $\frac{\partial \Sigma}{\partial u } =  -
 \Xi$, 
and integrating again: 
\be
x^A_{\ (B)} \ (t) \ = \ - \ (\frac{d_0}{r}) \ (\Sigma_{AB} \ (t) \ - \ \Sigma^-) \ . 
\ee
The limit $t \to \infty$ is taken and it follows that the test masses 
experience permanent displacements. 
Thus 
$\Sigma^+ - \Sigma^-$ is equivalent to an overall displacement of the test masses:
\be \label{reschrmem****4}
\triangle \ x^A_{\ (B)} \ = \ - \ (\frac{d_0}{r}) \ (\Sigma^+_{AB} \ - \
\Sigma^-_{AB}) \ . 
\ee
The right hand side of (\ref{reschrmem****4}) includes terms from the
electromagnetic field at highest order as given in our theorem \ref{displ*1}. 
Even though the form of (\ref{reschrmem****4}) is as in the EV case investigated by
Christodoulou in \cite{chrmemory} and \cite{chrdlbmathpgrt}, the nonlinear
contribution from the electromagnetic field is present in $\Sigma^+_{AB} -
\Sigma^-_{AB}$.

\newpage 

 
%
\vspace{10pt}
{\scshape Lydia Bieri \\ 
Department of Mathematics \\ 
University of Michigan \\ 
Ann Arbor, MI 48109, USA} \\ 
lbieri@umich.edu \\ \\ 
{\scshape PoNing Chen \\  
Department of Mathematics \\ 
Harvard University \\
Cambridge, MA 02138, USA} \\ 
pchen@math.harvard.edu \\ \\ 
{\scshape Shing-Tung Yau \\ 
Department of Mathematics \\ 
Harvard University \\
Cambridge, MA 02138, USA} \\ 
yau@math.harvard.edu

\end{document}